\documentclass[11pt]{amsart}
\usepackage{}
\usepackage{amssymb}
\theoremstyle{plain}
\newtheorem{Thm}{Theorem}

\begin{document}

\title[schrodinger equation and wave equation on finite graphs ]
{schrodinger equation and wave equation on finite graphs }

\author{Li Ma, Xiangyang Wang}
\address{L. Ma, Distinguished Professor, Department of mathematics \\
Henan Normal university \\
Xinxiang, 453007 \\
China} \email{nuslma@gmail.com}

\address{X.Y.Wang, Department of mathematics \\
Sun Yat Sen university \\
Guangzhou, 510275\\
China}

\email{mcswxy@mail.sysu.edu.cn}

\thanks{The research is partially supported by the National Natural Science
Foundation of China 10631020 and SRFDP 20090002110019}

\begin{abstract}
In this paper, we study the schrodinger equation and wave equation
with the Dirichlet boundary condition on a connected finite graph.
The explicit expressions for solutions are given and the energy
conservations are derived. Applications to the corresponding
nonlinear problems are indicated.


\end{abstract}

 \maketitle

\section{introduction}
In this paper we study the schrodinger equation and wave equation
with the Dirichlet boundary condition on a connected finite graph.
This work is a complement to the paper \cite{MW}. To our best
knowledge, this is the first paper on this direction.

A graph $G=(V,E)$ is a pair of the vertex-set $V$ and the edge-set
$E$. Each edge is an unordered pair of two vertices. If there is an
edge between $x$ and $y$, we write $x \sim y$. We assume that $G$ is
{\it local finite},
i.e., there exists a constant $c>0$ such that ${\rm deg}(x) := \#\{ y \in V: (x,y) \in E\} \le c$ for all $x\in G$. \\

Let $S$ be a finite subset of $V$, the {\it subgraph} $G(S)$
generated by $S$ is a graph, which consists of the vertex-set $S$
and all the edges $x \sim y, \ x,y \in S$ as the edge set.  The
boundary $\delta S$ of the induced subgraph $G(S)$ consists of all
vertices that are not in $S$ but adjacent to some vertex in $S$.We
assume that the subgraph $G(S)$ is connected. In below, we write
$(a,b)=a \overline{b}$ for complex numbers $a$ and $b$. We now
recall some facts from the book \cite{Ch}. Sometimes people may like
to write
$$
\bar{S}=S\bigcup \delta S.
$$

For a function $f: S\bigcup \delta S\to {\mathbb C}$, let
$$
\nabla_{xy} f=f(y)-f(x)
$$
for $y\sim x$. Then we recall that
$$
\int_{\bar{S}} f = \sum_{x \in \bar{S}} f(x), \quad
\|f\|^2=\int_{\bar{S}}|f|^2=\sum_{x\in \bar{S}} f(x)\overline{f(x)}
$$
and
$$
\|\nabla f\|^2=\sum_{\{x,y\in \bar{S}\}}|f(x)-f(y)|^2.
$$

We say that $f:S\bigcup \delta S\to {\mathbb C}$ satisfies the {\it Neumann boundary
condition} if for all $x\in \delta S$,
$$
\sum_{\{y\in S;y\sim x\}} (f(y)-f(x))=0.
$$

 Then the {\it Laplacian operator} can be written as
$$
(\Delta f)(x) = \sum_{y: y \sim x} (f(y) - f(x)) = \sum_{y: y \sim x} \nabla_{xy} f.
$$
Also the Neumann condition can be written as
$$
\sum_{\{y\in S;y\sim x\}} \nabla_{xy} f=0, \quad \forall x \in \delta S.
$$

We say that $f:S\bigcup \delta S\to {\mathbb C}$ satisfies the {\it Dirichlet
boundary condition} if  $f(x)=0$ for all $x\in \delta S$.

With Neumann or Dirichlet boundary condition, the Laplacian operator
on $S$ has finite eigenvalues  $\lambda_j>0$ with the corresponding
eigen-functions $\phi_j(x)$ \cite{Ch}, i.e.,
$$
-\Delta \phi_j=\lambda_j\phi_j, \ \ \ in \ \ S, \ \ \int_{\bar{S}}
|\phi_j|^2=1.
$$
In short, we can write this as
$$
-\Delta = \sum_j \lambda_jI_j
$$
where $I_j$ is the projection on to the j-th eigenfunction $\phi_j$
of the induced subgraph $S$ (see p.145 in \cite{Ch}).

Then we define the schrodinger kernel as
$$
\mathbf{S}_t(x,y)=\sum_jexp(-i\lambda_jt)\phi_j(x)\phi_j(y), \ \
t\geq 0.
$$
As in the heat kernel case we can write this as
$$
\mathbf{S}_t=\sum_j e^{-i\lambda_jt}I_j.
$$
Then
$$
\mathbf{S}_t=e^{it\Delta }=I+it\Delta+...., \ \ S_0=I.
$$

 For a function $f:S\bigcup \delta S\to {\mathbb C}$, we
define
\begin{equation}\label{sch}
u(x,t)=\sum_y\mathbf{S}_t(x,y)f(y), \ \  x\in S, \ t>0.
\end{equation}

Note that
$$
\mathbf{S}_0(x,y)=I=\sum_j\phi_j(x)\phi_j(y).
$$
Then for any $f:S\bigcup \delta S\to {\mathbb C}$,
$$
u(x,0)=\mathbf{S}_0(x,y)f(y)=\sum_y\sum_j\phi_j(x)\phi_j(y)f(y)=f(x).
$$
Then we can directly verify that the function $u$ satisfies the
schrodinger equation
$$
i\partial_tu(x,t)+\Delta u(x,t)=0, \ \ x\in S, \ t>0,
$$
with $u(x,0)=f(x)$ for $x\in S$. We denote $V(S)$ the space of
functions $f:S\bigcup \delta S\to {\mathbb C}$ satisfy the Dirichlet
boundary condition and with the $L^2$ norm.

 Then we show the following result.
\begin{Thm}\label{main} Assume that the function $f:S\bigcup \delta S\to {\mathbb C}$ satisfies the Dirichlet
boundary condition. Then there is a global solution $u(t):S\bigcup
\delta S\to {\mathbb C}$, which can be expressed in (\ref{sch}),  such that
$u$ satisfies the schrodinger equation
$$
i\partial_tu(x,t)+\Delta u(x,t)=0, \ \ x\in S, \ t>0,
$$
with $u(x,0)=f(x)$ for $x\in S$ and with the mass conservation
$$
\|u(t)\|^2=\|f\|^2
$$
and the energy conservation
$$
\int_{\bar{S}}|\nabla u(t)|^2=\int_{\bar{S}}|\nabla f|^2.
$$
\end{Thm}
The proof will be given in next section via the use of the spectrum
of the Laplacian and the Green formula.

Similarly we consider the wave equation on $\bar{S}$. Given two
functions $f$ and $g: \bar{S}\to R$ and both satisfy the Dirichlet
condition. We consider the following wave equation
\begin{equation}\label{wave}
u_{tt}=\Delta u(x,t), \ \ x\in S, \ t>0,
\end{equation}
with the initial conditions $u(x,0)=f(x)$ and $u_t(x,0)=g(x)$ for
$x\in S$.

We have the following result.
\begin{Thm}\label{wa} The solution to (\ref{wave}) is given by
\begin{equation}\label{solution} u(x,t)=\sum_j [I_jf(x)\cos
(\sqrt{\lambda_j} t)+I_jg(x)\frac{\cos (\sqrt{\lambda_j}
t)}{\sqrt{\lambda_j}}]\phi_j(x).
\end{equation}
Furthermore, we have  that
$$
\int_{\bar{S}}[|\nabla u|^2+u_t^2]dx=\int_{\bar{S}}[|\nabla
f|^2+g^2]dx.
$$
\end{Thm}

We remark that the related Duhamel principle can be used to give
solutions to the corresponding non-homogenous problems and we omit
them here. We don't consider problems on the infinite locally finite
connected graphs $V$ with symmetric weight as in \cite{MW} for the
following reason. One can prove the (skew) self-adjoint property of
$\Delta$ ($i\Delta$) (as proved in \cite{K}) which gives the
solutions to linear Schrodinger and wave equations on $L^2(V)$.
However, for such general graph $V$, the Strichartz type equality
for Schrodinger operator (or related interpolation inequality for
wave operator) is missing for applications to nonlinear Schrodinger
equations (for nonlinear wave equations).

The plan of this paper is below. In section \ref{sect2}, we prove
Theorem \ref{main}. We prove Theorem \ref{wa} in section
\ref{sect3}. We give some applications of our results in section
\ref{sect4}.

\section{Proof of Theorem \ref{main} }\label{sect2}

The importance of the boundary conditions above is the boundary term
vanishing in the formula below.

\begin{Thm}\label{green} Assume that $f:S\bigcup \delta S\to {\mathbb C}$. Then we have
\begin{equation}\label{green+1}
\int_S (\Delta f,f)=-\frac{1}{2}\int_{S}|\nabla f|^2+\sum_{x\in
S}\sum_{y\in\delta S}\overline{f(x)}\nabla_{xy}f.
\end{equation}
\end{Thm}
This can be verified directly. In fact, we can directly verify the
following more general formula (see Theorem 2.1 in \cite{G}) in a
compact form.

\begin{Thm}\label{green2} Assume that $f, \ g :S\bigcup \delta S\to {\mathbb C}$. Then we have
\begin{equation}\label{eq-thm-3}
\int_{\bar{S}} (\Delta f,g)=-\frac{1}{2}\sum_{x,y\in \bar{S}
}(\nabla_{xy} f, \nabla_{xy} g).
\end{equation}
\end{Thm}

We remark that the formula (\ref{eq-thm-3}) is proved in \cite{G}
for real functions, but the complex case can be done by writing the
complex function into the sum of real and imaginary parts.

For completeness, we give the proof of Theorem \ref{green2} below.
\begin{proof} As remarked above, we need only prove the result for
real functions.

 We make the following computation
$$
\begin{array}{ll}
\int_S (\Delta f, g) & = \sum_{x \in S} (\Delta f (x), g(x)) = \sum_{x \in S} \sum_{y: \ y \sim x} (\nabla_{xy} f, g(x))\\
& = \sum_{x \in S} \sum_{y\in S: \ y \sim x} (\nabla_{xy} f, g(x)) + \sum_{x \in S} \sum_{y\in \delta S: \ y \sim x} (\nabla_{xy} f, g(x))\\
& = \sum_{x, y \in S: x\sim y}(\nabla_{xy} f, g(y) - \nabla_{xy} g) + \sum_{x \in S} \sum_{y\in \delta S: \ y \sim x} (\nabla_{xy} f, g(x)) \\
& = - \sum_{x,y \in S, x \sim y} (\nabla_{xy} f, \nabla_{xy} g) + \sum_{x,y \in S, x \sim y} (\nabla_{xy} f, g(y)) \\
& + \sum_{x \in S} \sum_{y\in \delta S: \ y \sim x} (\nabla_{xy} f,
g(x)) =: A + B + C
\end{array}
$$
For the term $B$, we have
\[
\begin{array}{ll}
B & = \sum_{x,y \in S, x \sim y} (\nabla_{xy} f, g(y)) \\
& = \sum_{y \in S} (\sum_{x\in S \cup \delta S, x \sim y} \nabla_{xy} f, g(y)) - \sum_{y \in S} (\sum_{x\in \delta S, x \sim y} \nabla_{xy} f, g(y)) \\
& = \sum_{y\in S} (-\Delta f(y), g(y)) + \sum_{y \in S} (\sum_{x\in \delta S, x \sim y} \nabla_{yx} f, g(y)) \\
& = - \int_S (\Delta f, g) + C. \\
\end{array}
\]
(For the second term of the left hand side of the above equation, we
change $x$ to $y$,  $y$ to $x$, and see that this term is just $C$).
It follows that
\[
\int_S (\Delta f, g) = A + (- \int_S (\Delta f, g) + C) + C.
\]
Then we have
$$
2\int_S (\Delta f, g) = A  + C.
$$
We then re-write this into (\ref{eq-thm-3}). The proof is complete.
\end{proof}

For $f:S\bigcup \delta S\to {\mathbb C}$ satisfying the Dirichlet
condition, the boundary term in Theorem \ref{green} can be written
as
$$
-\sum_{x\in S}\sum_{y\in\delta
S}\overline{(f(y)-f(x))}\nabla_{xy}f=-\sum_{x\in S}\sum_{y\in\delta
S}|\nabla_{xy}f|^2,
$$
which is real. This fact is useful in the proof of mass and energy
conservation below.

 We now use Theorem \ref{green} (and Theorem
\ref{green2}) to prove Theorem \ref{main}.

\begin{proof}
Assume that $f$ satisfies either  Neumann or Dirichlet boundary
condition. Compute, via the use of the formula (\ref{green+1}),
$$
\frac{d}{dt}\|u(t)\|^2=2Re (u,u_t)=2Re(u, i\Delta u)=0.
$$
here we have used implicitly the boundary condition which implies
that $u_t=0$ on $\delta S$.
 We then have the mass conservation
$$
\|u(t)\|^2=\|f\|^2.
$$
Similarly we compute
$$ \frac{d}{dt}\int_{\bar{S}}|\nabla u(t)|^2=2Re (\nabla u,\nabla u_t).
$$
Using the formula (\ref{eq-thm-3}) for $f=u$ and $g=u_t$, we know
that
$$
Re (\nabla u,\nabla u_t)=-2Re\int_{\bar{S}} (\Delta u,u_t).
$$
Using the Schrodinger equation, we know that the term in the right
side is
$$ -2Re\int_{S} (\Delta u,i\Delta u),
$$
which is zero too.  Then
$$ \frac{d}{dt}\|\nabla u(t)\|^2=0
$$
 and we then get the energy conservation
$$
\|\nabla u(t)\|^2=\|\nabla f\|^2.
$$

The uniqueness of the solution follows from the mass conservation.
This completes the proof.
\end{proof}

\section{wave equations on finite graphs}\label{sect3}

Recall that two functions $f$ and $g: \bar{S}\to R$ satisfy the
Dirichlet condition.The problem under consideration is the following
wave equation
\begin{equation}\label{wave}
u_{tt}=\Delta u(x,t), \ \ x\in S, \ t>0,
\end{equation}
with the initial conditions $u(x,0)=f(x)$ and $u_t(x,0)=g(x)$ for
$x\in S$.

Here is the proof of Theorem \ref{wa}.
\begin{proof}
 Write the solution as
$u(x,t)=\sum_j u_j(t) \phi_j(x)$, where $u_j(t)=I_j u(x,t)$. Using
the initial conditions, we know that
$$
u_j(0)=I_jf(x), \ \ u_{jt}(0)=I_jg(x).
$$

Inserting $u(x,t)=\sum_j u_j(t) \phi_j(x)$ into (\ref{wave}) we get
that
$$
u_{jtt}+\lambda_ju_j=0.
$$
Solve this equation we know that
$$
u_j(t)=I_jf(x)\cos (\sqrt{\lambda_j} t)+I_jg(x)\frac{\cos
(\sqrt{\lambda_j} t)}{\sqrt{\lambda_j}}.
$$
Hence we can write the solution in the form (\ref{solution}).

Note that
$$
\frac{d}{dt}\int_{\bar{S}}[|\nabla
u|^2+u_t^2]dx=2\int_{\bar{S}}[(\nabla u, \nabla
u_t)+(u_{tt},u_t)]dx.
$$
Using the Green formula and (\ref{wave}) we then have
$$
\int_{\bar{S}}[(-\Delta u,u_t)+(u_{tt},u_t)]dx=0.
$$
That is,
$$
\frac{d}{dt}\int_{\bar{S}}[|\nabla u|^2+u_t^2]dx=0.
$$
This proves Theorem \ref{wa}.

\end{proof}

\section{discussions}\label{sect4}

  The advantage of our formulation of Schrodinger equation and wave equation
  in finite graphs is that it gives us the local existence result of
  the corresponding nonlinear equations by the use of fixed point
  theorem. For example, there is a unique local in time solution $u(x,t)$ to the nonlinear Schrodinger equation
\begin{equation}\label{non-sch}
i\partial_tu(x,t)+\Delta u(x,t)=|u|^{p-1}u(x,t), \ \ x\in S, \ t>0,
\end{equation}
with the initial data $u(x,0)=f(x)$ for $x\in S$ and with the
Dirichlet boundary condition. Here $p>1$. In fact, by Duhamel
principle, we know that the problem is equivalent to the fixed point
problem
$$
u(x,t)=\sum_y\mathbf{S}_t(x,y)f(y)+\int_0^t
\sum_y\mathbf{S}_{t-\tau}(x,y)|u|^{p-1}u(y,\tau)d\tau,
$$
on the Banach space $C^0([0,T],V(S))$. Here $T>0$ is sufficiently
small.

Using the Nehari method, one can easily obtain the following.
\begin{Thm}\label{Nehari} Given $p>1$ and $V(x)$ a non-negative function on $S$.
There exists a ground state solution $u\in V(S)$ to the problem
$$
-\Delta u(x)+V(x)u(x)=|u|^{p-1}u(x), \ \ u(x)>0, \ \ x\in S.
$$
\end{Thm}

Recall here that a ground state is a minimizer of the functional
$$
I(u)=-\frac{1}{2}\int_S (|\nabla
u(x)|^2+V(x)|u(x)|^2)-\frac{1}{p+1}\int_S |u(x)|^{p+1}
$$
over the set
$$
\{u\in V(S); u\neq 0,\int_S |u|^{p+1}=\int_S (|\nabla
u(x)|^2+V(x)|u(x)|^2)\}.
$$

 The minimizer exists and the proof is straightforward, so we omit it.

\bigskip

 {\bf Acknowledgement}. This work is done while both authors
are visiting the Department of Mathematics, CUHK, Hongkong in 2011
and the authors would like to thank the hospitality of the
Mathematical Department of CUHK.

\end{document}